\documentclass[11pt]{article}
\usepackage{amsmath}
\usepackage{amsfonts}
\usepackage{amssymb}
\newcommand{\ov}[1]{\overline{\vphantom{T}#1}}
\newcommand{\Z}{{\mathbb Z}}
\newcommand{\R}{{\mathbb R}}
\newcommand{\C}{{\mathbb C}}
\renewcommand{\th}{\theta}
\newcommand{\la}{\lambda}
\newcommand{\ga}{\gamma}
\newcommand{\eps}{\varepsilon}
\DeclareMathOperator{\cotan}{cotan}
\DeclareMathOperator{\sign}{sign}
\newcommand{\Proof}{{\it Proof. \/}}
\newcommand{\squareforqed}{\hbox{\rlap{$\sqcap$}$\sqcup$}}
\newcommand{\qed}{\ifmmode\squareforqed\else{\unskip\nobreak\hfil
\penalty50\hskip1em\null\nobreak\hfil\squareforqed
\parfillskip=0pt\finalhyphendemerits=0\endgraf}\fi}
\newcommand{\fp}{\qed\removelastskip\vskip\baselineskip\relax}
\newtheorem{theorem}{Theorem}[section]
\newtheorem{corollary}[theorem]{Corollary}
\newtheorem{proposition}[theorem]{Proposition}
\newtheorem{lemma}[theorem]{Lemma}
\newtheorem{definition}[theorem]{Definition}
\newtheorem{conjecture}[theorem]{Conjecture}

\begin{document}
\pagestyle{plain}
\title{Lambert $W$-Function Branch Identities}
\author{Henri Cohen}

\maketitle
\begin{abstract}
  After defining in detail the Lambert $W$-function branches, we give a large
  number of \emph{exact} identities involving (infinite) symmetric functions
  of these branches, as well as geometrically convergent series for all
  the branches. In doing so, we introduce a family of polynomials which
  may be of independent interest.
\end{abstract}

\bigskip

\section{Introduction: the Lambert $W$-Function Branches}

We begin by recalling the definition of the complex branches of the Lambert
$W$-function. Originally, the Lambert $W$-function was defined as the solution
of the implicit equation $we^w=z$ with $z\in\R$. It is immediate to see that
for $z\ge0$ there is a single solution, denoted $W_0(z)$, with a power series
expansion given by Lagrange inversion
$$W_0(z)=\sum_{n\ge1}(-1)^{n-1}\dfrac{n^{n-1}}{n!}z^n\;,$$
This power series converges for $|z|\le 1/e$ (here and in the sequel
$e=\exp(1)$), and provides an analytic continuation of $W_0(z)$ to the whole
interval $[-1/e,\infty[$. When $-1/e\le z<0$ there exists a \emph{second}
solution $W_{-1}(z)$ characterized by $W_{-1}(z)\le-1$. Note that
$W_0(-1/e)=W_{-1}(-1/e)=-1$ and that $W_{-1}(z)\to-\infty$ as $z\to0^-$.
Finally, if $z<-1/e$ there is no real solution.

\smallskip

For completeness and comparison with later formulas, note that Lagrange
inversion allows us to prove additional expansions such as
\begin{align*}W_0(z)^k&=\sum_{n\ge k}(-1)^{n-k}k\dfrac{n^{n-k-1}}{(n-k)!}z^n\text{\quad for all $k\in\Z\setminus\{0\}$\;,\ and}\\
\dfrac{1}{W_0(z)+1}&=\sum_{n\ge0}(-1)^n\dfrac{n^n}{n!}z^n\;.\end{align*}
More generally, if $F$ if any power series, we have
$$[z^n]F(W_0(z))=[t^n]F(t)(1+t)e^{-nt}\;,$$
where $[z^n]$ extracts the coefficient of $z^n$, see \cite{CGHJK}.

\smallskip

It is natural to extend these definitions to all complex numbers. Taking
logarithms in the initial implicit equation, we have $w+\log(w)-\log(z)=0$,
so for $k\in\Z$ we would like to denote by $w_k=W_k(z)$ the solution of the
equation $(E_k)$:
$$w_k+\log(w_k)=\log(z)+2k\pi i\;,$$
where $\log$ denotes the principal determination of the logarithm.
Note that this is \emph{not} the same as $w_k+\log(w_k/x)=2k\pi i$.

Because of the branching properties of the logarithm we will see that this
equation may have $0$, $1$, or $2$ solutions, so we must be more careful.
Although the definitions are in the literature, for instance in
\cite{CGHJK}, we prefer doing everything from scratch. The precise result is
as follows:\footnote{We evidently use the
  much more sensible French notation $]a,b[$ and similar for intervals, instead
  of the dreadful $(a,b)$ which has so many other different meanings.}

\begin{proposition} Let $z\in\C\setminus\{0\}$, $k\in\Z$, and denote
  by $(E_k)$ the equation $w+\log(w)=\log(z)+2k\pi i$.
  \begin{enumerate}
  \item If $z\notin]-\infty,0[$, equation $(E_k)$ has a unique complex
    solution denoted $W_k(z)$, and we have $W_{-k}(z)=\ov{W_k(\bar{z})}$.
  \item If $z\in]-\infty,0[$ and $k\ne0,-1$, or $z\in]-\infty,-1/e[$
    and $k=0$ or $k=-1$, equation $(E_k)$ has a unique complex solution
    denoted $W_k(z)$, and we have $W_{-k}(z)=\ov{W_{k-1}(z)}$.
  \item If $z\in[-1/e,0[$ and $k=-1$, equation $(E_k)$ has no solution
    (this is the only case in which there is no solution).
  \item If $z\in[-1/e,0[$ and $k=0$, equation $(E_k)$ has two real solutions
    denoted $W_{-1}(z)$ and $W_0(z)$ characterized by
    $W_{-1}(z)\le-1\le W_0(z)$, which coincide if $z=-1/e$.
\end{enumerate}
\end{proposition}

We give the detailed elementary proof, since it is sometimes obscured
in the literature. It can also be proved using the principle of the
argument, but the proof is essentially the same, see \cite{Khe}.

\Proof As usual denoting by $\log(y)$ the principal branch of the logarithm,
we can write $\log(z)=a+\la \pi i$ with $a\in\R$ and $-1<\la\le1$. Equation
$(E_k)$ is $w+\log(w)=a+(2k+\la)\pi i$, so writing $w=\rho e^{i\th}$ with
$-\pi<\th\le\pi$ gives the system $S_{2k+\la}(a)$, where $S_{\mu}(a)$ is
$$\rho\sin(\th)+\th=\mu\pi\text{\quad and\quad}\rho\cos(\th)+\log(\rho)=a\;.$$

Let us first find the solutions with $\th\in]0,\pi[$. 
The first equation has the unique solution $\rho=(\mu\pi-\th)/\sin(\th)$,
and since we have $\rho>0$, we must have $\th<\mu\pi$. Thus if $\mu\le0$
there are no solutions with $\th\in]0,\pi[$, so assume $\mu>0$.
Replacing in the second equation gives
$$F_{\mu}(\th):=(\mu\pi-\th)\cotan(\th)+\log(\mu\pi-\th)-\log(\sin(\th))=a\;.$$
We compute that
$$F'_{\mu}(\th)=-\dfrac{\mu\pi-\th-\sin(2\th)}{\sin^2(\th)}-\dfrac{1}{\mu\pi-\th}\;,$$
and a tedious computation shows that $F'_{\mu}(\th)<0$ for all
$\th\in]0,\pi[$ such that $\th<\mu\pi$. It follows that $F_{\mu}$ is a
strictly decreasing function of $\th$.

When $\th\to0^+$ we have $F_{\mu}(\th)\to+\infty$. If $\mu>1$, when
$\th\to\pi^-$ we have $F_{\mu}(\th)\to-\infty$, so $F_{\mu}$ has
a single root, and it follows that $S_{\mu}(a)$ has a single solution.
If $\mu<1$, when $\th\to\mu\pi^-$ we again have $F_{\mu}(\th)\to-\infty$,
so again $S_{\mu}(a)$ has a single solution. Finally, if $\mu=1$, when
$\th\to\pi^-$ we compute that $F_1(\th)\to-1$. It follows that
if $a>-1$, $S_1(a)$ has a single solution, while if $a\le-1$ it has no
solution.

To summarize, $S_{\mu}(a)$ has a solution with $\th\in]0,\pi[$ if and only if
$\mu>0$ and either $\mu\ne1$ or $\mu=1$ and $a>-1$, and that solution is
unique.

Since changing $\th$ into $-\th$ in the system simply changes $\mu$ into
$-\mu$, we deduce that $S_{\mu}(a)$ has a solution with $\th\in]-\pi,0[$
if and only if $\mu<0$ and either $\mu\ne-1$ or $\mu=-1$ and $a>-1$,
and that solution is unique.
    
A solution of $S_{\mu}(a)$ with $\th=0$ is possible if and only if $\mu=0$ and
$\rho+\log(\rho)=a$, and since this is an increasing function of $\rho$ and
tends to $-\infty$ when $\rho\to0^+$ and to $+\infty$ when $\rho\to\infty$, we
deduce that there always exists a solution with $\th=0$ when $\mu=0$, and it is
unique (and evidently equal to $W_0(e^a)$).

A solution of $S_{\mu}(a)$ with $\th=\pi$ is possible if and only if $\mu=1$
and $-\rho+\log(\rho)=a$. This is a concave function of $\rho$ tending to
$-\infty$ when $\rho\to0^+$ and $\rho\to\infty$, with a maximum at $\rho=1$
equal to $-1$, so there exists a solution with $\th=\pi$ if and only if
$\mu=1$ and $a\le-1$, and there are two solutions if $a<-1$ (evidently
equal to $W_0(-e^a)$ and $W_{-1}(-e^a)$), and a double solution if $a=-1$.

\smallskip

Summarizing, we have proved the following:

\begin{lemma} The solutions to $S_{\mu}(a)$ are as follows:
  \begin{enumerate}
  \item If $\mu>1$, or $0<\mu<1$, a unique solution with $\th\in]0,\pi[$,
    if $\mu<-1$ or $-1<\mu<0$ a unique solution with $\th\in]-\pi,0[$,
    and in both cases no solution with $\th=0$ or $\th=\pi$.
  \item If $\mu=-1$ and $a>-1$, a unique solution with $\th\in]-\pi,0[$,
    and if $\mu=-1$ and $a\le-1$ no solution.
  \item If $\mu=1$ and $a>-1$, a unique solution with $\th\in]0,\pi[$,
    and if $\mu=1$ and $a\le-1$ two solutions with $\th=\pi$ which coincide for $a=-1$.
  \item If $\mu=0$, a unique solution with $\th=0$.
\end{enumerate}\end{lemma}

Since $\mu=2k+\la$ with $-1<\la\le1$ and $k\in\Z$, $\mu>1$ is equivalent to
$k\ge1$, $\mu<-1$ is equivalent to $k\le -1$, $-1<\mu<1$ and $\mu\ne0$
is equivalent to $k=0$ and $\la\ne0$, $\mu=0$ is equivalent to $k=\la=0$,
$\mu=1$ is equivalent to $k=0$ and $\la=1$, and $\mu=-1$ is equivalent to
$k=-1$ and $\la=1$, so it is clear that this lemma is equivalent to the
proposition.\fp

\begin{corollary} Keep the same notation.
  \begin{enumerate}
  \item If $z\notin]-\infty,0[$, or $z\in]-\infty,0[$ and $k\ne0,-1$,
    or $z\in]-\infty,-1/e[$ and $k=0$ or $k=-1$, $W_k(z)$
    is the unique solution of $(E_k)$, and $W_{-k}(z)=\ov{W_k(\bar{z})}$
    if $z\notin]-\infty,0[$, and $W_{-k}(z)=\ov{W_{k-1}(z)}$ otherwise.
  \item If $z\in[-1/e,0[$, $W_{-1}(z)$ is the unique real solution of $E_0$
    \emph{(not of $E_{-1}$)}
      less than or equal to $-1$, and $W_0(z)$ is the unique real solution
      of $E_0$ greater or equal to $-1$.\end{enumerate}\end{corollary}
  
To give a precise idea of what these results mean, let us give the
example of $W_k(x)$ for $-2\le k\le 1$, and $x$ in each of the three real
intervals $]-\infty,-1/e[$, $[-1/e,0[$, and $]0,\infty[$, for instance
$x=-1$, $x=-0.1$, and $x=1$:
\begin{align*}
  W_{-2}(-1)&=-2.06227... - 7.58863...\quad W_{-2}(-0.1)=-4.44909... - 7.30706...\\
  W_{-1}(-1)&=-0.31813... - 1.33723...\quad W_{-1}(-0.1)=-3.57715...\\
  W_0(-1)&=-0.31813... + 1.33723...\quad\ \, W_0(-0.1)=-0.11183...\\  
  W_1(-1)&=-2.06227... + 7.58863...\quad\ \, W_1(-0.1)=-4.44909... + 7.30706...
\end{align*}
In both of these cases, $W_{-k}(z)=\ov{W_{k-1}(z)}$, except in the second
case for $k=0$ and $1$.
\begin{align*}
  W_{-2}(1)&=-2.40158... - 10.77629...\\
  W_{-1}(1)&=-1.53391... - 4.37518...\\
  W_0(1)&=\phantom{-}0.56714...\\ 
  W_1(1)&=-1.53391... + 4.37518...
\end{align*}
In this case, $W_{-k}(z)=\ov{W_k(z)}$.

\medskip

{\bf Remark.} It is possible to define $W_k(z)$ for \emph{any} $k\in\C$
(not only for $k\in\Z$) by the defining equation $W_k(z)+\log(W_k(z))=\log(z)+2k\pi i$.
Of course, if $k\notin\Z$, this will not anymore correspond to a solution
of $ye^y=z$. Furthermore, it is clear that we can find $\th$ such that
$-1/2<\th\le 1/2$ and such that
$W_k(z)+\log(W_k(z))=\log(ze^{2\th\pi i})+2m\pi i$,
where $m\in\Z$ is such that $|m-k|<1$, so if for instance
$z\notin]-\infty,0]$, the solution to this equation will simply be
$W_m(ze^{2\th\pi i})$. A similar remark applies to the equation
$w+\log(w/x)=2k\pi i$ whose solution is $W_{k+\eps}(x)$ for some $\eps=-1$,
$0$, or $1$ depending on $k$ and $x$.

\section{The Fundamental Identity}

We begin by the following lemma:

\begin{lemma} Let $x\in\C\setminus\{0\}$ be fixed. We have
  $$\lim_{t\to\pm\infty}\sum_{k\in\Z}\dfrac{1}{W_k(x)-t}=\mp\dfrac{1}{2}\;,$$
  where here and below a sum or product for $k\in\Z$ is understood as the
  limit as $K\to\infty$ of the sum or product for $-K\le k\le K$.
\end{lemma}

\Proof We will see below the trivial fact that $W_k(x)=2k\pi i+O(\log(k))$
when $|k|\to\infty$,
so the series indeed converges if we sum symmetrically. Denoting it by $f(t)$,
we thus have $f(t)=f_1(t)+f_2(t)$ with
$$f_1(t)=\sum_{k\in\Z}\dfrac{1}{2k\pi i-t}\text{\quad and\quad}f_2(t)=\sum_{k\in\Z}\dfrac{2k\pi i-W_k(x)}{(W_k(x)-t)(2k\pi i+t)}\;.$$
Since $2k\pi i-W_k(x)=O(\log(k))$ the series for $f_2(t)$ converges uniformly
in $t$, hence $\lim_{t\to\pm\infty}f_2(t)=0$. On the other hand it is classical
that $f_1(t)=-\coth(t/2)/2$, so the result follows.\fp

\begin{theorem}
  For $x\in\C\setminus\{0\}$ we have
  $$\prod_{k\in\Z}\left(1-\dfrac{t}{W_k(x)}\right)=e^{-t/2}-\dfrac{t}{x}e^{t/2}\;.$$
\end{theorem}

\Proof The variable $x$ being fixed, consider the function 
$F(t)=e^{-t/2}-(t/x)e^{t/2}$. It has order $1$, has no pole, and its zeros
are the roots of $te^t=x$, so by definition they are the $W_k(x)$. By the
Hadamard factorization theorem it thus has a Hadamard product which we can
write in the form
$$F(t)=Ae^{Bt}\prod_{k\in\Z}(1-t/W_k(x))\;.$$
Note that the standard way of writing the product would be
$\prod_{k\in\Z}(1-t/W_k(x))e^{t/W_k(x)}$, but since we know that the
product written as above converges if we interpret it as the limit
of the product for $|k|\le K$, it is preferable to write it as we have done.

We trivially have $A=F(0)=1$. To determine $B$, we compute the logarithmic
derivative:
$$\dfrac{F'(t)}{F(t)}=B-\sum_{k\in\Z}\dfrac{1}{W_k(x)-t}\;,$$
We immediately compute that $\lim_{t\to\pm\infty}F'(t)/F(t)=\pm1/2$,
so making $t\to\infty$ or $t\to-\infty$ and using the lemma gives $B=0$,
proving the theorem.\fp

{\bf Remarks.}\begin{enumerate}\item
In the proofs of the lemma and the theorem, we have implicitly
used the evident fact that for $t$ real with $|t|$ sufficiently large,
for all $k\in\Z$ we have $t\ne W_k(x)$ and $t\ne 2k\pi i$.
\item Since $W_k(x)$ is close to $2k\pi i$, the above theorem should be
compared with the identity (equivalent to the one we used above for
$\coth(t/2)$):
$$-t\prod_{k\in\Z\setminus\{0\}}\left(1-\dfrac{t}{2k\pi i}\right)=e^{-t/2}-e^{t/2}\;.$$
\end{enumerate}

The following corollary gives the usual form of the Hadamard product:

\begin{corollary} We have
  $$\prod_{k\in\Z}\left(1-\dfrac{t}{W_k(x)}\right)e^{t/W_k(x)}=e^{t/x}\left(1-\dfrac{t}{x}e^t\right)\;.$$
\end{corollary}

\Proof Immediate from the theorem since identification of the coefficients
of $t^1$ gives $\sum_{k\in\Z}1/W_k(x)=1/2+1/x$.\fp
  
\smallskip

Using the theorem, differentiating with respect to $t$ or $x$, and/or
specializing, we easily find a large number of identities which we regroup
in the following corollary:

\begin{corollary} We have the following identities:
$$\sum_{k\in\Z}\dfrac{1}{W_k(x)-t}=\dfrac{1}{2}+\dfrac{1+t}{xe^{-t}-t}\;,\quad \sum_{k\in\Z}\dfrac{1}{W_k(x)+1}=\dfrac{1}{2}\;,$$
$$\sum_{k\in\Z}\dfrac{1}{W_k(x)}=\dfrac{1}{2}+\dfrac{1}{x}\;,\quad
\sum_{k\in\Z}\dfrac{1}{W_k(x)^2}=\dfrac{2x+1}{x^2}\;,\quad
\sum_{k\in\Z}\dfrac{1}{W_k(x)^3}=\dfrac{3x^2+6x+2}{2x^3}\;.$$
$$\prod_{k\in\Z}\left(1-\dfrac{t}{W_k(x)}\right)=e^{-t/2}-\dfrac{t}{x}e^{t/2}\;,\quad
\sum_{k\in\Z}\dfrac{1}{(W_k(x)+1)(W_k(x)-t)}=\dfrac{1}{xe^{-t}-t}\;.$$
$$\sum_{k\in\Z}\dfrac{1}{(W_k(x)-t)^2}=\dfrac{xe^{-t}(2+t)+1}{(xe^{-t}-t)^2}\;,\quad\sum_{k\in\Z}\dfrac{1}{(W_k(x)+1)^2}=\dfrac{1}{xe+1}\;.$$
$$\sum_{k\in\Z}\dfrac{1}{(W_k(x)+1)(W_k(x)-t)^2}=\dfrac{xe^{-t}+1}{(xe^{-t}-t)^2}\;,\quad\sum_{k\in\Z}\dfrac{1}{(W_k(x)+1)^3}=\dfrac{1}{xe+1}\;.$$
$$\sum_{k\in\Z}\dfrac{W'_k(x)}{W_k(x)-t}=\dfrac{1}{x}\left(\dfrac{1}{2}+\dfrac{t}{xe^{-t}-t}\right)\;,\quad\sum_{k\in\Z}\dfrac{W'_k(x)}{W_k(x)}=\dfrac{1}{2x}\;.$$
\end{corollary}

It is clear that if $F$ is any rational function, by decomposing into partial
fractions, any symmetrically convergent series of the form
$\sum_{k\in\Z}F(W_k(x))$ can be evaluated explicitly. For instance,
\begin{align*}\sum_{k\in\Z}\dfrac{1}{W_k(x)^2+t^2}&=\dfrac{x(\cos(t)+\sin(t)/t)+1}{x^2+2xt\sin(t)+t^2}\text{\quad and}\\
  \sum_{k\in\Z}\dfrac{W_k(x)}{W_k(x)^2+t^2}&=\dfrac{1}{2}+\dfrac{x(\cos(t)-t\sin(t))-t^2}{x^2+2xt\sin(t)+t^2}\;.\end{align*}

\bigskip

{\bf Acknowledgment:} I am grateful to I.~Pinelis for having given the
identity $\sum_{k\in\Z}1/(W_k(x)+1)=1/2$ in the {\tt mathoverflow} forum,
which led to the present section of this paper, but his proof is completely
different, using an integral representation of this sum coming from the
residue theorem.

\section{Asymptotics and Sum and Product Formulas}

In this section we are going to show that even nonconvergent symmetric
sums or products of the $W_k$ can be evaluated explicitly after removing
the asymptotic main term. First note the following easy proposition:

\begin{proposition}\label{prop:asymp} The variable $x$ being fixed, as
$|k|\to\infty$ we have
$$W_k(x)=2k\pi i-\log(2k\pi i)+\log(x)+\dfrac{\log(2k\pi i)-\log(x)}{2k\pi i}+O\left(\dfrac{\log(|k|)^2}{k^2}\right)\;.$$
\end{proposition}

\Proof For notational simplicity set $K=2k\pi i$, $L=\log(2k\pi i)=\log(K)$,
$X=\log(x)$, and $M=L-X=\log(2k\pi i)-\log(x)$. By successive approximation, we
find that if we set
$w_k=K-M+M/K+\eps=K(1-M/K+M/K^2+\eps/K)$ we have
$\log(w_k)=L-M/K+O(L^2/K^2)$, hence
$w_k+\log(w_k)=K+X+\eps+O(L^2/K^2)$, and using that
$w_k+\log(w_k)=\log(x)+2k\pi i$ we see that $\eps=O(L^2/K^2)$.\fp

Note that we will give below the complete asymptotic expansion.

\begin{corollary}\label{cor:asymp} Let the variable $x$ be fixed.
  \begin{enumerate}
  \item For $k>0$ sufficiently large we have
  $(2k-1)\pi<|W_k(x)|<(2k+1)\pi$.
  \item If $x$ is positive real, then for $K$ sufficiently large, the solutions
    to the equation $we^w=x$ with $|w|<(2K+1)\pi$ are the $W_k(x)$ for
    $-K\le k\le K$.\end{enumerate}
\end{corollary}

\Proof Immediate from the proposition and left to the reader.\fp

Note that if $x$ is not a positive real, there may be an additional solution
for $k=\pm(K+1)$.
  
\begin{theorem}\label{thm:cons} We have
  $$\lim_{K\to\infty}\left(\sum_{-K\le k\le K}W_k(x)+\log(\pi^{2K}(2K)!)-(2K+1/2)\log(x)\right)=\dfrac{\log(2)}{2}\;.$$
\end{theorem}

\Proof We will prove this theorem in two stages. First we will prove that
the left-hand side is a constant independent of $x$. Second, the
appendix due to Letong Hong and Shengtong Zhang proves that this constant
is equal to $\log(2)/2$, which I conjectured in the first version of this
paper. Note that their proof also shows that the left-hand side is independent
of $x$, but we have chosen to give the present proof first since it handles
more easily the case where $x$ is not a positive real.

\smallskip

For $k>0$ we have
$$W_k(x)+W_{-k}(x)=-2\log(2k\pi)+2\log(x)+1/(2k)+O(\log(k)^2/k^2)\;.$$
It follows that
$$\sum_{-K\le k\le K}W_k(x)=-2\log(K!)+2K\log(x/(2\pi))+\log(K)/2+F(x)+o(1)$$
for some function $F$. By Stirling's formula, we have
\begin{align*}-2\log(K!)-&2K\log(2\pi)+\log(K)/2=-(2K+1/2)\log(2K)+2K\\
  &\phantom{=}-2K\log(\pi)-\log(4\pi)+o(1)\\
  &=-\log((2K)!)-2K\log(\pi)-\log(4\pi)+\log(2\pi)/2+o(1)\;,\end{align*}
so
$$\sum_{-K\le k\le K}W_k(x)=-\log((\pi/x)^{2K}(2K)!)+A(x)+o(1)$$
for some other function $A(x)$. To find $A(x)$, we note that\footnote{This is not quite rigorous, but one can prove that the derivative of the $o(1)$ is still
  a $o(1)$.}
$$A'(x)+\dfrac{2K}{x}+o(1)=\sum_{-K\le k\le K}W'_k(x)\;.$$
Now since $W'_k(x)=(1/x)W_k(x)/(W_k(x)+1)=(1/x)(1-1/(W_k(x)+1))$, we have
by Pinelis's formula
$$\sum_{-K\le k\le K}W'_k(x)=\dfrac{2K+1/2}{x}+o(1)\;,$$
so $A'(x)=1/(2x)$, hence $A(x)=\log(x)/2+A$ for some constant $A$.
The fact that $A$ is real follows by choosing $x>0$ and using that
$W_{-k}(x)=\ov{W_k(x)}$, proving that the left-hand side is indeed a real
constant. We refer to the appendix for the proof that $A=\log(2)/2$ and
a second proof of the fact that the left-hand side is constant.\fp

The proof shows that the implicit $o(1)$ term is a $O(\log(K)^2/K)$.

\begin{corollary} We have
  $$\lim_{K\to\infty}\dfrac{\prod_{-K\le k\le K}W_k(x)}{\pi^{2K}(2K)!}=\sqrt{x/2}\;,$$
    and more generally
$$\lim_{K\to\infty}\dfrac{\prod_{-K\le k\le K}(W_k(x)-t)}{\pi^{2K}(2K)!}=\sqrt{x/2}(e^{-t/2}-(t/x)e^{t/2})\;,$$    
\end{corollary}

\Proof By definition, we have $\log(W_k(x))=\log(x)+2k\pi i-W_k(x)$. It follows
from the theorem that
\begin{align*}\sum_{-K\le k\le K}\log(W_k(x))&=(2K+1)\log(x)-\sum_{-K\le k\le K}W_k(x)\\
  &=\log(\pi^{2K}(2K)!)+\log(x)/2-\log(2)/2+o(1)\;,\end{align*}
so the result follows for $t=0$, and the more general result from the
formula for $\prod_k(1-t/W_k(x))$ given above.\fp

\section{A Family of Polynomials}

Before going further, we are now going to study in great detail a family of
polynomials important for our work. These polynomials date back at least
to Comtet \cite{Com}.

\begin{definition} Let $A(X)$ be a given polynomial such that $A(0)=0$.
  \begin{enumerate}\item We define the Lambert family of polynomials
    $P_n(X)=P_{n,A}(X)$ associated to $A$ by $P_0(X)=A(X)$ and for $n\ge0$
    by the recursion
    $$P_{n+1}(X)=-P_n(X)+n\int_0^XP_n(t)\,dt\;.$$
  \item We denote by $F_A(T,X)=\sum_{n\ge0}P_{n,A}(X)T^n$ the generating
    function of the $P_n(X)$.\end{enumerate}
\end{definition}

\newcommand{\pd}{\partial}

\begin{lemma}\label{lem0} $P_n$ is the Lambert family associated to $A$ if
  and only if
  $$(1+T)\dfrac{\pd F_A}{\pd X}=A'(X)+T^2\dfrac{\pd F_A}{\pd T}$$
  and $P_n(0)=0$ for all $n$.\end{lemma}

\Proof Since $\pd F_A/\pd X=\sum_{n\ge0}P'_n(X)T^n$ we have
\begin{align*}(1+T)\pd F_A/\pd X&=dP_0/dX+\sum_{n\ge1}(P'_n(X)+P'_{n-1}(X))T^n\\
  &=A'(X)+\sum_{n\ge1}(n-1)P_{n-1}(X)T^n=A'(X)+T^2\pd F_A/\pd T\;.\end{align*}
Conversely, if this is satisfied we have $P'_0(X)=A'(X)$ and
$P'_n(X)+P'_{n-1}(X)=(n-1)P_{n-1}(X)$, and since we assume in addition that
$P_n(0)=0$ we have $P_0(X)=A(X)$ and $P_n(X)+P_{n-1}(X)=(n-1)\int_0^XP_{n-1}(t)\,dt$.\fp

\begin{lemma}\label{lem1} If $P_n$ is the Lambert family associated to $A(X)$
then $\int_0^XP_n(t)\,dt$ is the Lambert family associated to
$\int_0^XA(t)\,dt$.
\end{lemma}

\Proof If $Q_n(X)=\int_0^XP_n(t)\,dt$ it is clear that
$Q'_{n+1}(X)+Q'_n(X)=nQ_n(X)$, so the result follows since $Q_n(0)=0$.\fp

\begin{theorem}\label{thpn} Let $L_n$ be the Lambert family associated to
  $A(X)=-X$, $F(T,X)=\sum_{n\ge0}L_n(X)T^n$ its generating function, and
  $k\ge1$. By abuse of notation, write $F^k(T,X)=\sum_{n\ge0}L_{n,k}(X)T^n$
  (including for $k\le-1$), so that $L_{n,1}(X)=L_n(X)$.
  \begin{enumerate}
  \item For $n\ge1$ we have
    $$L_{n,2}(X)=-2(n+1)\int_0^XL_n(t)\,dt=-\dfrac{2(n+1)}{n}(L_{n+1}(X)+L_n(X))\;.$$
      Equivalently, $L_{n,2}(X)/(n+1)$ is the Lambert family associated
      to\\ $A(X)=X^2$.
    \item The function $F$ satisfies the differential equation in $T$
      $$TFF'+(1+T)F'+F=0\;.$$
    \item The function $F$ satisfies the implicit equation
      $$\log(1+TF)+F=-X\;.$$
    \end{enumerate}
\end{theorem}

\Proof (1). By Lemma \ref{lem0} applied to $A(X)=-X$ we have
\begin{align*}(1+T)\pd(F^k)/\pd X&=kF^{k-1}(1+T)\pd F/\pd X\\&=kF^{k-1}(A'(X)+T^2\pd F/\pd T)=-kF^{k-1}+T^2\pd(F^k)/\pd T\;.\end{align*}
Identifying coefficients of $T^n$, it follows that
$$L'_{n,k}(X)+L'_{n-1,k}(X)=(n-1)L_{n-1,k}(X)-kL_{n,k-1}(X)\;.$$
We prove (1) by induction. Since $L_{0,2}(X)=(-X)^2=(-2)\int_0^X(-t)\,dt$,
it is true for $n=0$. Assume that it is true for $n-1$: we thus have
$$L'_{n,2}(X)=2nL_{n-1,1}(X)-2(n-1)n\int_0^XL_{n-1,1}(t)\,dt-2L_{n,1}(X)$$
and since $(n-1)\int_0^XL_{n-1,1}(t)\,dt=L_{n,1}(X)+L_{n-1,1}(X)$ we obtain
$L'_{n,2}(X)=-2(n+1)L_{n,1}(X)$, proving our induction hypothesis since
$L_{n,k}(0)=0$ for all $n$, $k$. It follows from Lemma \ref{lem1} that
$L_{n,2}(X)/(n+1)$ is the Lambert family associated to $X^2$.

\smallskip

(2). For $n\ge1$ the coefficient of $T^n$ in
$TFF'+(1+T)F'+F=(T/2)(F^2)'+F'+TF'+F$ is equal to
\begin{align*}&\phantom{=}(n/2)L_{n,2}(X)+(n+1)L_{n+1,1}+nL_{n,1}+L_{n,1}
\\&=-n(n+1)\int_0^XL_{n,1}(t)\,dt+(n+1)n\int_0^XL_{n,1}(t)\,dt=0\;,\end{align*}
and the coefficient of $T^0$ is $L_{0,1}(X)+L_{1,1}(X)=-X+X=0$, proving
the differential equation.

\smallskip

(3). If we set $G=\log(1+TF)+F$, we have
\begin{align*}G'&=(TF'+F)/(1+TF)+F'\\
  &=((1+T)F'+F+TFF')/(1+TF)=0\end{align*}
by (2), so $G(T)=G(0)=F(0)=-X$.\fp

{\bf Remarks.}\begin{enumerate}
\item The proof of (3) shows that the general solution to the differential
  equation in (2) is a solution to the implicit equation
  $\log(1+TF)+F=A(X)$, where $A(X)=F(0)$.
\item It is clear that the theorem can be proved by using
    Lagrange inversion, but I have preferred to give the above direct proof.
  \end{enumerate}

It is natural to call the $L_n(X)$ the \emph{Lambert polynomials}.
  
\begin{proposition}\label{prop:pde} In addition to the partial differential
  equation $(1+T)\pd F/\pd X=-1+T^2\pd F/\pd T$ of Lemma \ref{lem0}, the
  function $F$ also satisfies the following ones:
  $$\dfrac{\pd F}{\pd X}=-\dfrac{TF+1}{TF+T+1}\text{\quad and\quad}F\dfrac{\pd F}{\pd X}=(TF+1)\dfrac{\pd F}{\pd T}=-T\dfrac{\pd F}{\pd T}-F\;.$$
  \end{proposition}

\Proof By Lemma \ref{lem0} and Theorem \ref{thpn} (2) we have
\begin{align*}\dfrac{\pd F}{\pd X}&=\dfrac{1}{1+T}\left(-1+T^2\dfrac{\pd F}{\pd T}\right)\\
&=-\dfrac{1}{1+T}\left(1+T^2\dfrac{F}{TF+T+1}\right)=-\dfrac{TF+1}{TF+T+1}\;.\end{align*}
Thus, by Theorem \ref{thpn} once again we have $TF+T+1=-F/F'$,
hence $\pd F/\pd X=(TF+1)(\pd F/\pd T)/F$, and the last formula
again follows from the theorem.\fp

\begin{corollary}\label{cor:log} We have the identity
  $$\sum_{n\ge1}\dfrac{L'_n(X)}{n}T^n=-\log\left(\dfrac{F(T,X)}{-X}\right)\;.$$
\end{corollary}

\Proof Indeed, if we differentiate with respect to $T$, the derivative of
the left hand side is $\sum_{n\ge1}L'_n(X)T^{n-1}=(\pd F/\pd X+1)/T$,
while that of the right-hand side is $-(\pd F/\pd T)/F$, and equality follows
from the last formula given by the proposition. The identity follows since
both sides vanish for $T=0$,\fp

\begin{corollary}\begin{enumerate}
    \item Generalizing (1) of the theorem, for $k\ge2$ we have
      \begin{align*}L_{n,k}(X)&=-\dfrac{k}{k-1}(n+k-1)\int_0^XL_{n,k-1}(t)\,dt\\
        &=-\dfrac{k}{k-1}\dfrac{(n+1)L_{n+1,k-1}(X)+(n+k-1)L_{n,k-1}(X)}{n}\;,\end{align*}
      and $L_{n,k}(X)/\binom{n+k-1}{n}$ is the Lambert family associated to
      $A(X)=(-X)^k$.
    \item For all $k\ge2$ we have
      \begin{align*}L_{n,k}(X)&=(-1)^{k-1}k(k-1)\binom{n+k-1}{n}\int_0^X(X-t)^{k-2}L_n(t)\,dt\\
        &=(-1)^{k-1}k\binom{n+k-1}{n}\int_0^X(X-t)^{k-1}L'_n(t)\,dt\;.\end{align*}
    \item For any polynomial $A(X)$ such that $A(0)=0$, the Lambert family
      associated to $A$ is given by
      \begin{align*}P_{n,A}(X)&=-\int_0^XA'(X-t)L'_n(t)\,dt\\&=-A'(0)L_n(X)-\int_0^XA''(X-t)L_n(t)\,dt\;.\end{align*}
\end{enumerate}\end{corollary}

\Proof (1). Multiplying by $F^{k-2}$ the second identity of the proposition
we obtain $F^{k-1}\pd F/\pd X=-TF^{k-2}\pd F/\pd T-F^{k-1}$, and identifying
the coefficients of $T^n$ gives $L'_{n,k}/k=-((n+k-1)/(k-1))L_{n,k-1}$,
proving the first formula since $L_{n,k}(0)=0$.

It is clear that $(n!/(n+k-1)!)L_{n,k}(X)$ is a Lambert family for $k=1$.
Assume that this is the case for $k-1$ with $k\ge2$. Using the formula just
proved and using the induction hypothesis, we have
\begin{align*}&\phantom{=}\dfrac{n!}{(n+k-1)!}L'_{n,k}(X)+\dfrac{(n-1)!}{(n+k-2)!}L'_{n-1,k}(X)\\
&=-\dfrac{k}{k-1}\left(\dfrac{n!}{(n+k-2)!}L_{n,k-1}(X)+\dfrac{(n-1)!}{(n+k-3)!}L_{n-1,k-1}(X)\right)\\&=-\dfrac{k}{k-1}(n-1)\dfrac{(n-1)!}{(n+k-3)!}\int_0^XL_{n-1,k-1}(t)\,dt\\&=(n-1)\dfrac{(n-1)!}{(n+k-2)!}L_{n-1,k}(X)\;,\end{align*}
proving that $(n!/(n+k-1)!)L_{n,k}(X)$ is a Lambert family, and it is
associated to $L_{0,k}(X)/(k-1)!=(-X)^k/(k-1)!$. Using this and replacing
in the first formula proves the second, proving (1).

\smallskip

(2). This follows immediately from the formula giving iterated integrals.

\smallskip

(3). By linearity, it is clear that
$$P_{n,A}(X)=-A'(0)L_n(X)-\int_0^XA''(X-t)L_n(t)\,dt\;,$$
so the result follows by integrating by parts.\fp

The polynomials $L_{n,k}$ can be given explicitly:
\begin{proposition}\label{prop:stir} We have
$$L_{n,k}(X)=(-1)^{n+k}k!\binom{n+k-1}{n}\sum_{j=1}^ns(n,n+1-j)\dfrac{X^{k-1+j}}{(k-1+j)!}\;,$$
where the $s(n,m)$ are the Stirling numbers of the first kind. In particular,
\begin{align*}L_n(X)&=(-1)^{n+1}\sum_{j=1}^ns(n,n+1-j)\dfrac{X^j}{j!}\;,\text{\quad or equivalently,}\\
\int_0^\infty e^{-Xt}L_n(t)\,dt&=X^{-(n+1)}\prod_{j=1}^{n-1}(j-X)\;.\end{align*}
\end{proposition}

\Proof The result for $k=1$ follows immediately from the recursion for
Stirling numbers $s(n+1,k)=s(n,k-1)-ns(n,k)$, the general result by
successive integration, and the integral formula from the definition of
Stirling numbers.\fp

{\bf Remark.} In view of the above results, it would seem that the most
natural polynomials are the $L'_n(X)$ and not the $L_n(X)$: they are monic,
have constant term equal to $(-1)^{n-1}$, satisfy the same Lambert-type
recursion (but they are not a Lambert family since $L'_n(0)\ne0$), but it is
immediate to see that the equations involving them are more complicated, so
we have preferred to take the $L_n(X)$ themselves as fundamental building
blocks. We will see in the next section that the higher
derivatives of $L_n$ enter in the expansion of $F^k$ for $k\le-1$.

\smallskip

The first few $L_n(X)$ are the following:
\begin{align*}L_0(X)&=-X\;,\ L_1(X)=X\;,\ L_2(X)=X^2/2-X\;,\\
  L_3(X)&=X^3/3-(3/2)X^2+X\;,\ L_4(X)=X^4/4-(11/6)X^3+3X^2-X\;.\end{align*}

\section{Further Results on Lambert Polynomials}

The preceding section dealt essentially with the explicit computation of
the polynomials $L_{n,k}$ when $k\ge1$. The present section deals with
the case $k\le-1$.

\begin{proposition} We have the additional partial differential equation
  $$\dfrac{\pd^2 F}{\pd X^2}=-T^2\dfrac{\pd^2(1/F)}{\pd T^2}\;.$$
\end{proposition}

\Proof By Proposition \ref{prop:pde} we have
$\pd F/\pd X=-(TF+1)/(TF+T+1)$. Differentiating once more we deduce that
\begin{align*}\dfrac{\pd^2 F}{\pd^2 X}&=-\dfrac{T^2}{(TF+T+1)^2}\dfrac{\pd F}{\pd X}=T^2\dfrac{TF+1}{(TF+T+1)^3}\;.\end{align*}

On the other hand, writing for simplicity ${}'$ instead of
$\pd/{\pd T}$, by Theorem \ref{thpn} (2) we have
$(1/F)'=-F'/F^2=1/(F(TF+T+1))$. Now,
\begin{align*}(TF^2+TF+F)'&=F^2+2TFF'+F+TF'+F'=F^2+TFF'\\
  &=F^2-TF^2/(TF+T+1)=F^2(TF+1)/(TF+T+1)\;,\end{align*}
hence $(1/F)''=-(TF+1)/(TF+T+1)^3$, proving the proposition.\fp

\begin{corollary}\label{cor:lneg} Recall that we have defined $L_{n,k}$ for
  all $k\in\Z$ by $F^k(T,X)=\sum_{n\ge0}L_{n,k}(X)T^n$.
  \begin{enumerate}\item
    We have $L_{0,-1}(X)=L_{1,-1}(X)=-1/X$, and for $n\ge2$
    $$L_{n,-1}(X)=-\dfrac{L_n''(X)}{n(n-1)}\;,\text{\quad in other words}$$
    $$\dfrac{1}{F(T,X)}=-1/X-(1/X)T-\sum_{n\ge2}\dfrac{L_n''(X)}{n(n-1)}T^n\;.$$
  \item More generally, let $k\le-2$. We have
    $$(n+k)L_{n,k}(X)=-\dfrac{k}{(k+1)}L'_{n,k+1}(X)\;,$$
    and for $n\ge 1-k$ we have
    $$L_{n,k}(X)=(-1)^{k-1}k\dfrac{(n+k-1)!}{n!}L^{(1-k)}_n(X)\;.$$
  \end{enumerate}
\end{corollary}

\Proof (1). The formulas for $L_{0,-1}(X)$ and $L_{1,-1}(X)$ are obtained
directly, and the formula for $L_{n,-1}(X)$ is equivalent to the PDE of
the proposition. It can also be easily obtained by differentiating the
formula of Corollary \ref{cor:log}.

\smallskip

(2). Multiplying by $F^{k-1}$ the last PDE of Proposition \ref{prop:pde}, we
obtain $F^k\pd F/\pd X=-TF^{k-1}\pd F/\pd T-F^k$, and identifying the coefficients
of $T^n$ gives $L'_{n,k+1}/(k+1)=-((n+k)/k)L_{n,k}$, proving the first
formula, and the second follows by induction.\fp

The above corollary does not give $L_{n,k}$ for $n\le -k$. For this,
we need to introduce another family of polynomials, closely linked
to $L_n$ as follows:

\begin{definition} For $n\ge1$ we define a family $M_n$ of functions
  by $M_1(X)=-1/X$ and $M_{n+1}(X)=M'_n(X)/n-M_n(X)$.
\end{definition}

Note that the $M_n$ are polynomials in $1/X$ (for instance $M_2(X)=(X+1)/X^2$
and $M_3(X)=-(2X^2+3X+2)/(2X^3)$), but the derivative $M'_n(X)$ is of course
taken with respect to $X$, not with respect to $1/X$.

\begin{proposition} Let $k\le-1$.
  \begin{enumerate} \item We have
  $$L_{n,k}=\begin{cases}
    -\dfrac{k}{n}\dfrac{M_n^{(-k-n)}(X)}{(-k-n)!}\text{\quad for $1\le n\le -k$\;,}\\
    (-X)^k\text{\quad for $n=0$.}\end{cases}$$
\item
  In other words, for $k\le-1$ we have
  \begin{align*}F(T,X)^k&=(-X)^k-k\sum_{1\le n\le -k}\dfrac{M_n^{(-k-n)}(X)}{n(-k-n)!}T^n\\&\phantom{=}+(-1)^{k-1}k\sum_{n\ge1-k}\dfrac{(n+k-1)!}{n!}L^{(1-k)}_n(X)T^n\;.\end{align*}
  \end{enumerate}
\end{proposition}
  
\Proof In the proof of Theorem \ref{thpn} (1) we have seen the recursion
$$L'_{n,k}(X)+L'_{n-1,k}(X)=(n-1)L_{n-1,k}(X)-kL_{n,k-1}(X)\;,$$
which is valid for all $k$ including negative ones. Thus, for $n\ge1$ we have
$$L'_{n+1,-n}(X)+L'_{n,-n}(X)=nL_{n,-n}(X)+nL_{n+1,-n-1}(X)\;.$$
On the other hand, the recursion given in (2) of the above corollary implies
that $L'_{n+1,-n}(X)=0$ for $n\ge1$. We thus have the recursion
$L_{n+1,-n-1}(X)=L'_{n,-n}(X)/n-L_{n,-n}(X)$, which is the recursion
for $M_n$, proving the result for $n=-k$ since $L_{1,-1}(X)=-1/X=M_1(X)$.

The same recursion for $1\le n<-k$ gives
$L_{n,k}(X)=(-k/((k+1)(k+n)))L'_{n,k+1}(X)$, so by
induction for $1\le j\le -k-n$ we have, since $n\ge1$ hence $j\le -k-1$:
$$L_{n,k}(X)=\dfrac{k}{k+j}\dfrac{(-k-j-n)!}{(-k-n)!}L^{(j)}_{n,j+k}\;,$$
and choosing $j=-k-n$ and using what we just proved gives the result for
$1\le n\le -k$. The case $n=0$ is trivial directly since $L_{0,k}$ is the
coefficient of $T^0$ in $F^k$, hence equal to $(-X)^k$.\fp

The functions $M_n$ can also be given explicitly in terms of Stirling numbers:

\begin{proposition}
$$M_n(X)=\dfrac{1}{(n-1)!}\sum_{j=0}^n(-1)^{j-1}\dfrac{s(n,j+1)j!}{X^{j+1}}
    =\dfrac{(-1)^n}{(n-1)!}\int_0^\infty e^{-Xt}\prod_{j=1}^{n-1}(t+j)\,dt\;.$$
\end{proposition}

\Proof Left to the reader.\fp

The coefficients of $M_n(X)$ are of course closely related to those of
$L_n(X)$, but I do not see any natural way of expressing $M_n(X)$ in
terms of the $L_n(X)$, except by the following formula coming from the
definition of the beta function:

\begin{proposition}
  $$M_n(X)=\dfrac{n(n+1)}{X^{n+1}}\int_0^\infty\dfrac{L_n(-Xt)}{(t+1)^{n+2}}\,dt\;.$$
\end{proposition}

\Proof Simply use the expansions of $L_n(X)$ and $M_n(X)$ in terms of
Stirling numbers. In addition, one can integrate once or twice by parts
and obtain similar formulas involving $L'_n$ and $L''_n$.\fp

\section{Convergent Series for the Lambert Branches}

The result of Proposition \ref{prop:asymp} is asymptotic, and can be
pushed further if desired. The result is that not only does it give an
asymptotic expansion to any desired number of terms, but in fact a
\emph{convergent series} for $W_k(x)$ even for small $k$. We can even
put these series in a slightly more general setting as follows.

\begin{proposition}\label{prop:asympk} Fix $x\in\C\setminus\{0\}$, let
  $K\in\C\setminus\{0\}$, and define
  $$M=K+\log(K)-\log(x)-2k\pi i\;.$$
\begin{enumerate}\item As $|K|\to\infty$ we have the asymptotic expansions
  $$W_k(x)=K+\sum_{n\ge0}\dfrac{L_n(M)}{K^n}\text{\quad and\quad}\log(W_k(x))=\log(K)-\sum_{n\ge1}\dfrac{L_n(M)}{K^n}\;,$$
  where the $L_n(z)$ are as usual the Lambert polynomials.
\item These series converge to $W_k(x)$ and $\log(W_k(x))$ when $2(1+|M|)<|K|$
  (we will see below the exact domain of convergence in special cases).
\item More generally, we have the following formulas:
  \begin{align*}(W_k(x)-K)^j&=\sum_{n\ge0}\dfrac{L_{n,j}(M)}{K^n}\text{\quad for all $j\in\Z\setminus\{0\}$\;,}\\
    \dfrac{-M}{W_k(x)-K}&=1+\dfrac{1}{K}+\sum_{n\ge2}\dfrac{ML''_n(M)}{n(n-1)K^n}\;,\\
    \dfrac{1}{1+W_k(x)}&=\sum_{n\ge1}\dfrac{L'_n(M)}{K^n}\;,\text{ and }
    \log\left(\dfrac{W_k(x)-K}{-M}\right)=-\sum_{n\ge1}\dfrac{L'_n(M)}{nK^n}\end{align*}
  with the same convergence properties.\end{enumerate}
\end{proposition}

\Proof (1). Denote by $R$ the right-hand side of the expansion for $W_k(x)$,
and recall that $F(T,X)=\sum_{n\ge0}L_n(X)T^n$. By Theorem \ref{thpn}, we have
$$\log(R)=\log(K)+\log(1+(1/K)F(1/K,M))=\log(K)-M-F(1/K,M)\;.$$
It follows that
$$R+\log(R)=K+\log(K)-M=\log(x)+2k\pi i\;,$$
which is exactly equation $(E_k)$ defining the Lambert branches, and since
$W_k(x)$ is the solution of $(E_k)$ when $k\ne0,1$, we deduce that it gives
its asymptotic expansion. The expansion for $\log(W_k(x))$ follows from
$W_k(x)+\log(W_k(x))=\log(x)+K$.

\smallskip

(2). Let us study the convergence properties. By the recursion formula for
Stirling numbers, it is immediate to show that $|s(n,k)|\le 2^nn!/k!$, so
replacing in the explicit formula for $L_n(X)$ we deduce that
$|L_n(M)|\le 2^n(1+|M|)^n$, proving the convergence result.
Using the stronger bound $|s(n,n-k)|\le(n/2)^k\binom{n-1}{k}$
conjectured by Ren\'e Gy and proved by Mike Earnest in the {\tt stackexchange}
forum, we can replace $2^n$ by $e^{n/2}$.

\smallskip

(3). These formulas follow from (1), the definition, and using
$W'_k(x)=W_k(x)/(x(1+W_k(x)))$ and $L'_0(X)=-1$.\fp

\begin{corollary}\begin{enumerate}\item
  We have the asymptotic expansions
  $$W_k(x)=K+\sum_{n\ge0}\dfrac{L_n(M)}{K^n}$$
  (as well as all the others given in the proposition)
  for instance for $K=2k\pi i$ and $M=\log(K)-\log(x)$
  (which gives an asymptotic expansion when $|k|\to\infty$),
  or for $K=2k\pi i+\log(x)$ and $M=\log(K)$
  (which gives another asymptotic expansion when $|k|\to\infty$, and also
  when $|x|\to\infty$ or $x\to0$).
\item If $x\in[-1/e,0[$, we have the same expansions also for $k=-1$
    by choosing $K=\log(-x)$ and $M=\log(-K)$.
\item All these expansions converge if and only if
  $$\log(|K|)+1-\Re(M)+\min_{m\in\{-1,0\}}\Re(W_m(-e^{M-1}))>0$$
  (note that this condition is automatically satisfied in case (2)).
\end{enumerate}\end{corollary}

\Proof (1) is clear, (2) is proved in exactly the same way as the
proposition, and (3) is proved in \cite{CGHJK} and \cite{Kal-Jef}.\fp

{\bf Remarks.}\begin{enumerate}
\item The series with $K=2k\pi i+\log(x)$ is the one usually given,
for instance in \cite{Com} and \cite{CGHJK}. However, when $k\ne0$ the series
with $K=2k\pi i$ often converges faster.
\item We could choose $K=W_k(x)$, in which case $M=0$ by definition,
so the identity is a triviality since $L_n(0)=0$ for all $n\ge0$.
\item Even for $k=1$ (which is certainly not
in the ``asymptotic'' regime) and $x=1$, we can compute $W_1(1)$ to $38$
decimals using slightly more than $100$ terms of the series (of course the
usual iterative methods to compute $W_k(x)$ are much more efficient since they
are based on Newton or Halley iterations).
\item We have remarked above that it is possible, although not very useful,
to define $W_k(x)$ even when $k\notin\Z$. It is clear that by construction
the above series are still valid in this more general context.
\end{enumerate}

\medskip

In \cite{Kal-Jef}, some alternate expansions are given for the principal
branch $W_0$, but they are trivially generalizable to all branches. For
instance, instead of choosing $K$ as main variable, we can choose $K_1=K+1$.
Thus, we write
\begin{align*}W_k(x)&=K+F(1/K,M)=K_1-1+F(1/(K_1-1),M)\\
  &=K_1-1+\sum_{n\ge0}L_n(M)/(K_1-1)^n\\
  &=K_1-1+\sum_{n\ge0}L_n(M)/K_1^n\sum_{j\ge0}\binom{n+j-1}{n-1}/K_1^j\\
  &=K_1-1+\sum_{N\ge0}(1/K_1^N)\sum_{n=0}^N\binom{N-1}{n-1}L_n(M)\;,\end{align*}
so
$$W_k(x)=K_1-1+\sum_{N\ge0}\dfrac{P_n(M)}{K_1^N}\text{\quad with\quad}P_n(X)=\sum_{n=0}^N\binom{N-1}{n-1}L_n(X)\;,$$
and these new polynomials $P_n(X)$ can be expressed using another type of
Stirling numbers, and also have recursion properties coming from those of
$L_n(X)$, such as
$$P_{n+1}(X)=\int_0^X(nP_n(t)-(n-1)P_{n-1}(t))\,dt\;,$$
whose proof is left to the reader. The first few $P_n(X)$ are
\begin{align*}P_0(X)&=-X\;,\ P_1(X)=X\;,\ P_2(X)=X^2/2\;,\ P_3(X)=X^3/3-X^2/2\;,\\
P_4(X)&=X^4/4-(5/6)X^3\;,\ P_5(X)=X^5/5-(13/12)X^4+X^3/2\;.\end{align*}
In particular, note that $P_N(X)$ is divisible by $X^{\lfloor N/2\rfloor + 1}$.

For the branch $k=0$, \cite{Kal-Jef} mentions that the domain of convergence
is much larger than that of the original series, but for the other branches
the difference is not that large since $K+1$ is not so different from $K$.

\section{Additional Properties and Conjectures}

\begin{conjecture}\begin{enumerate}
    \item The roots of $L_n(X)$ are all nonnegative real, and less
      than or equal to $n+\log(n)+\ga+O(1/n)$, where $\ga$ is Euler's
      constant.
    \item More generally, the $j$th root of $L_n(X)$ in decreasing order is
      asymptotic to $n/j+\sum_{j\le m\le n}1/m$.
    \item The roots of $L_n(X)$ and $L_{n+1}(X)$ interlace.
    \item More generally (and more precisely), for $k\ge1$ the roots of
      $L_{n,k}(X)$ are all nonnegative real, and their maximum has an
      asymptotic expansion
      $$n+\log(n)+\ga+k-1+C_0(k)/n+C_1(k)/n^2+...\;,$$
      where $C_m(k)$ is a polynomial of degree $m$ in $k$, in particular
      $C_0(k)$ is independent of $k$. This is still true for
      $k\le0$ if we interpret $L_{n,k}(X)$ as $L_n^{(1-k)}(X)$.
      \end{enumerate}
\end{conjecture}

Note that the fact that the real roots of $L_{n,k}(X)$ are nonnegative
follows immediately from Proposition \ref{prop:stir} and the fact that
Stirling numbers alternate in sign, and the fact that the roots of
$L_{n+1}(X)$ and $L_n(X)$ interlace follows from the recursion once proved that
all the roots of $L_n(X)$ are real.

Numerically,

$$C_0(k)=-4.3469909700078207218721533281574751800...$$

The above conjecture is certainly intimately linked to the formula
$$\int_0^\infty e^{-Xt}L_n(t)\,dt=X^{-(n+1)}\prod_{j=1}^{n-1}(j-X)$$
given by Proposition \ref{prop:stir}, but I do not yet see a proof.

\section{Integral Representation of the Branches}

The following is part of a result given in \cite{Khe} and is given for
completeness:

\begin{proposition}
  \begin{enumerate}\item
    Set $P=\pi/2+1$, $K=2k\pi i+\log(x)$, and $f(t)=t-\log(t)+K$.
    Assume that either $k\ne0,-1$, or that $k=0$ or $k=-1$ and
    $x\notin[-1/e,0]$. We have $W_k(x)=N/D$, with
    \begin{align*}
      N&=\dfrac{K}{K^2+P^2}+\int_0^\infty\dfrac{tdt}{(t^2+1)(f(t)^2+\pi^2)}\\
      D&=\dfrac{P}{K^2+P^2}-\int_0^\infty\dfrac{dt}{(t^2+1)(f(t)^2+\pi^2)}\;.\end{align*}
  \item Set $P=\pi/2-1$, $K=(2k-\sign(k))\pi i+\log(x)$, and $f(t)=t+\log(t)-K$.
    Assume that $k\ne-1,0,1$. Then $W_k(x)=-N/D$ with the same formulas as
    in (1).\end{enumerate}
\end{proposition}

{\bf Remarks}\begin{enumerate}
\item I refer to \cite{Khe} for an explicit but much more complicated formula
  for $k=0$ or $k=-1$ and $x\in[-1/e,0[$.
\item The conditions given by the author for the validity of these
  integral formulas are not quite the same as those given here.
  I believe at least that the above are correct.
\item As the author mentions, the change of variable $t=\sinh(u)$
      transforms the integrals into exponentially decaying ones, which can
      thus be easily computed. But in fact the doubly-exponential integration
      method does this automatically.
\end{enumerate}

\vfill\eject

\centerline{\large \bf Appendix: Proof of Theorem \ref{thm:cons}}

\bigskip

\newcommand{\abs}[1]{\left\vert #1 \right \vert}

\centerline{By Letong Hong and Shengtong Zhang}

\bigskip

Without loss of generality assume that $x$ is a fixed positive real number\footnote{This is purely for simplicity. The argument works for any $x\in \C \backslash \{0\}$ with minor modifications.}, and set $f(w)=we^w-x$. By
Corollary \ref{cor:asymp}, for $K$ sufficiently large the set of zeros $z$ of
$f$ with $\abs{z}\le (2K+1)\pi$ is precisely $\{w_k:-K \leq k \leq K\}$,
where we set $w_k=W_k(x)$. By Jensen's formula applied to $f$ we have
\begin{align*}\log x&=\sum_{-K\le k\le K}\log\abs{w_k}-(2K+1)\log ((2K+1)\pi)+\frac{I((2K+1)\pi)}{2\pi}\;,\text{ with}\\
I(\rho)&=\int_0^{2\pi}\log\abs{f(\rho e^{i\theta})}\,d\theta\;.\end{align*}
Since $w_ke^{w_k} = x$, we have $\log \abs{w_k} = \log x - \Re w_k$, so
summing over $k\in [-K,K]$ we obtain
\begin{equation}\label{eq1}
\sum_{-K\le k\le K} \Re w_k = 2K\log x - (2K+1)\log ((2K+1)\pi) + \frac{I((2K+1)\pi)}{2\pi}\;.\end{equation}
We now study the integral $I((2K+1)\pi)$ when $K$ is large, so we set
$\rho(K)=(2K+1)\pi$. The function $f(w) = we^w - x$ behaves like $we^w$ when
$\Re w \geq 0$ and is near $-x$ when $\Re w \leq -2 \log \rho(K)$. We now handle
the intermediate case.
\begin{lemma}
For $K$ sufficiently large and any $w$ with $\abs{w} =\rho(K)$, we have
$$\abs{f(w)} \gg 1\;.$$
\end{lemma}

\Proof 
If $\abs{we^w} \notin [x/2, 2x]$ the result follows by the triangle inequality.
Otherwise, the real part of $w$ must be $O(\log x + \log(\rho(K)))$, therefore
the imaginary part of $w$ must be $\pm\rho(K) + o(1)$. Since $\rho(K)$ is an odd
multiple of $\pi$, it follows that the argument of $we^w$ must be
$\pm \frac{\pi}{2} + o(1)$. Since $x$ is positive real, this implies that
$\abs{we^w - x}$ is also bounded from below.\fp

We can now give a precise estimate for $I((2K+1)\pi)$:
\begin{lemma}
For $K$ sufficiently large, we have
\begin{align}
  \int_0^{2\pi}\log\abs{f(\rho(K) e^{i\theta})} d\theta = \pi \log \rho(K) + 2\rho(K) + \pi \log  x + o(1) \;.
\end{align}
\end{lemma}

\Proof For simplicity, write $\rho$ instead of $\rho(K)$, set $w = \rho e^{i\theta}$, and let $\epsilon = \frac{3\log \rho}{\rho}$. We split the integral into three terms corresponding to $\theta \in [0, \frac{\pi}{2}] \cup [\frac{3\pi}{2}, 2\pi]$, $\theta\in [\frac{\pi}{2}, \frac{\pi}{2} + \epsilon] \cup [\frac{3\pi}{2} - \epsilon, \frac{3\pi}{2}]$, and $\theta \in [\frac{\pi}{2} + \epsilon, \frac{3\pi}{2} - \epsilon]$. In the first term, we have $\Re w \geq 0$, so
$\abs{we^w} = \rho e^{\Re w} \geq \rho$. It follows that
$\log\abs{f(w)} = \log \abs{we^w} + O(\rho^{-1})$, hence
\begin{align*}
 \int_{\theta \in [0, \frac{\pi}{2}] \cup [\frac{3\pi}{2}, \pi]}\log\abs{f(\rho e^{i\theta})} &= \int_{\theta \in [0, \frac{\pi}{2}] \cup [\frac{3\pi}{2}, \pi]}\log \abs{we^w} + O(\rho^{-1})    \\
 &=  \int_{\theta \in [0, \frac{\pi}{2}] \cup [\frac{3\pi}{2}, \pi] }(\log \rho + \Re w) + O(\rho^{-1}) \\
 &= \pi \log \rho + 2\rho + O(\rho^{-1}).
\end{align*}
In the second term, since $w$ has negative real part, we have the upper bound
$\abs{f(w)} \leq  x + \abs{we^w} \leq 2\rho$, and by the previous lemma,
it is bounded from below by a constant. It follows that $\log \abs{f(w)} \ll \log \rho$,
so the second term tends to $0$.

Finally, in the third term we have $\Re w \leq -2\log \rho$, so that
$\abs{we^w} = \rho e^{\Re w} \ll \rho^{-1}$, therefore we have
$\log \abs{f(w)} = \log  x + O(\rho^{-1})$. It follows that
$$\int_{\theta \in [\frac{\pi}{2} + \epsilon, \frac{3\pi}{2} - \epsilon]}\log\abs{f(\rho e^{i\theta})} = \pi \log  x + o(1).$$
Summing the three terms proves the lemma.\fp

Plugging this estimate into \eqref{eq1}, and using that for $x>0$ we have
$\Im(w_{-k})=-\Im(w_k)$, we deduce that as $K\to \infty$:
$$
\sum_{-K\le k\le K}W_k(x) =(2K+\frac12)\log  x - (2K+\frac12)\log ((2K+1)\pi) + (2K+1) + o(1).$$
Using Stirling's formula, it is immediate to deduce that
$$
\sum_{-K\le k\le K}W_k(x) =(2K+\frac12)\log  x -\log(\pi^{2K}(2K)!)+\dfrac{\log(2)}{2}+o(1)\;,$$
proving Theorem \ref{thm:cons}.

\end{document}